\documentclass[12pt]{article}

\usepackage{amsmath,amsthm,amssymb,amsfonts}
\usepackage{latexsym}
\def\pmb#1{\setbox0=\hbox{$#1$}%
\kern-.025em\copy0\kern-\wd0
\kern.05em\copy0\kern-\wd0
\kern-.025em\raise.0433em\box0}

\newcommand{\oM}{\overline{M}}
\newcommand{\ra}{\rightarrow}

\newcommand{\la}[1]{\mbox{\large $#1$}}
\newcommand{\La}[1]{\mbox{\Large $#1$}}
\newcommand{\LA}[1]{\mbox{\LARGE $#1$}}

\newcommand{\h}{\mbox{\lie h}}                  
\newcommand{\R}{\mbox{\field \symbol{82}}}          

\newfont{\lie}{eufm10 at 12pt}
\newfont{\field}{msbm10 at 11pt}

\newtheorem{theorem}{Theorem}[section]
\newtheorem{lemma}{Lemma}[section]
\newtheorem{corollary}{Corollary}[section]
\newtheorem{proposition}{Proposition}[section]

\begin{document}
\title{\Large \bf Heinz mean curvature estimates  in warped product
spaces $M\times_{e^{\psi}}N$}
\author{ Isabel M.C.\ Salavessa}
\date{}
\protect\footnotetext{\!\!\!\!\!\!\!\!\!\!\!\!\! {\bf MSC 2010:}
Primary: 53C21, 53C40; 53C24; Secondary: 58C40,
\\
Supported by the FCT through the project UID/CTM/04540/2017.}
\maketitle ~~~\\[-5mm]
\begin{quotation}\noindent
{\footnotesize   Center of Physics and Engineering of Advanced Materials (CeFEMA), Instituto 
Superior T\'{e}cnico, University
of Lisbon, Edif\'{\i}cio Ci\^{e}ncia, Piso 3, Av.\ Rovisco Pais,
1049-001 Lisboa, Portugal;~~
e-mail: isabel.salavessa@ist.utl.pt}\\[5mm]
\hspace*{8cm} \em In grateful memory of James Eells \em \\[2mm]
{\small {\bf Abstract:}  If a graph submanifold $(x,f(x))$  of a Riemannian warped 
product space  $(M^m\times_{e^{\psi}}N^n,\tilde{g}=g+ e^{2\psi}h)$ is immersed with 
parallel mean curvature $H$, then we obtain a Heinz type estimation of the mean curvature. Namely, on each compact domain $D$ of $M$,
 $m\|H\|\leq \frac{A_{\psi}(\partial D)}{V_{\psi}(D)}$ holds, 
where $A_{\psi}(\partial D)$ and $V_{\psi}(D)$ are the ${\psi}$-weighted
area and volume, respectively. 
In particular, $H=0$  if $(M,g)$ has zero weighted Cheeger constant, a concept
recently introduced by D.\ Impera et al.\ (\cite{[Im]}). This generalizes the known cases
 $n=1$ or $\psi=0$. We also conclude minimality using a closed calibration, 
assuming  $(M,g_*)$ is complete where $g_*=g+e^{2\psi}f^*h$,
and for some constants $\alpha\geq \delta\geq 0$, $C_1>0$ and  $\beta\in [0,1)$,
 $\|\nabla^*\psi\|^2_{g_*}\leq \delta$, $\mathrm{Ricci}_{\psi,g_*}\geq \alpha$, 
and $det_g(g_*)\leq C_1 r^{2\beta}$  holds when $r\to +\infty$,  where $r(x)$ is the distance  function on $(M,g_*)$ from some fixed point.
Both results rely on expressing the squared norm of the mean curvature as a weighted divergence of a suitable vector field.}
\end{quotation}

\markright{\sl\hfill  Salavessa \hfill}

\section{Introduction}
\renewcommand{\thesection}{\arabic{section}}
\renewcommand{\theequation}{\thesection.\arabic{equation}}
\setcounter{equation}{0}
In 1955, E.\ Heinz  \cite{[H]}  obtained an estimative of the 
mean curvature of a piece of surface of $\R^3$ described by a graph of a function in terms 
of an isoperimetric inequality. Using some integral formulas, he proved that, if $f(x,y)$ 
is a function defined on the disc $x^2+ y^2<R^2$ and the mean curvature of the surface 
$(x,y,f(x,y))$ satisfies $\|H\|\geq c>0$, where $c$ is a constant, then $R\leq \frac{1}{c}$. Thus, if $f$ is defined in all $\R^2$ and  $\|H\|$ is constant, then $H=0$. Ten years later 
this problem was extended and solved for the case of a function $f:\R^m \ra \R$, by Chern 
\cite{[Ch]}, and independently by Flanders \cite{[F]}.  
 In 1985, James Eells proposed us to find out if such isoperimetric estimations on the mean 
curvature could be extended  in higher codimensions.
Given a map  $f:M\ra  N$ between two Riemannian manifolds $(M,g)$ and $(N,h)$, of dimensions 
$m$ and $n$, respectively, assuming the graph submanifold  $\Gamma_f:=\{(x,f(x)):x\in M \}$ 
of the Riemannian product $\oM=(M\times N, g\times h)$ has  parallel mean curvature $H$, we proved 
the following inequality  
(\cite{[S1]}, \cite{[S2]})
\begin{equation}\label{old}
 \|H\|\leq \frac{1}{m}\frac{A(\partial D)}{V(D)},
\end{equation}
holds on each  compact domain $D\subset M$, 
where $A(\partial D)$ and  $V(D)$ are respectively the area of $\partial D$ and the volume of $D$ with respect to the metric $g$. Since $\|H\|$ is constant, we have 
\begin{equation}
\label{Cheeger}
\|H\|\leq \frac{1}{m}\h(M,g),
\end{equation}
where
$$\label{Cheeger2}
\h(M,g)=\inf_D\frac{A(\partial D)}{V(D)},
$$
is the Cheeger constant of $(M,g)$, with $D$ ranging  over all open domains of $M$ with compact closure in $M$ and smooth boundary (see e.g.\ \cite{[Cha]}). This constant is zero, if, for example, $M$ is a closed manifold or it is a simple noncompact Riemannian manifold, i.e., there exists a diffeomorphism $\phi:(M,g)\ra (\R^{m}, <,>)$ onto $\R^m$
such that $\lambda g \leq \phi^*<,>\leq \mu g$ for some positive constants $\lambda, \mu$.
Another large class of spaces with zero Cheeger constant are the 
complete Riemannian manifolds with nonnegative Ricci tensor (cf.\ \cite{[ADR]}). 
Inequality (\ref{Cheeger}) is  sharp. In \cite{[S1],[S2],[S3], [LiSa]} it is given 
examples of graphic hypersurfaces in $\mathbb{H}^m\times \mathbb{R}$, with constant mean curvature $c$ any real in $[0,m-1]$, where $\mathbb{H}^m$ is the hyperbolic space,  a space
with Cheeger constant equal to $(m-1)$. These examples are of the form $(x, f(r(x)))$ where
$r(x)$ is the distance function to a point of $\mathbb{H}^m$.

This problem was generalized by Guanghan Li and the author in \cite{[LiSa]}, in the  
context of calibrated manifolds. Given a  $m$-form defining a calibration $\Omega$, not necessarily closed,  on a Riemannian manifold $(\oM,\bar{g})$ of dimension $m+n$,
and  $F:M\ra \oM$  an oriented immersed submanifold of dimension $m$, the 
\em $\Omega$-angle \em of  $M$ is the function $\cos\theta:M\to [-1,1]$ given by
\begin{equation}\label{angle}
\cos\theta=\Omega(X_1,\ldots, X_m),
\end{equation}
where $X_i$ is a direct o.n. frame of $T_pM$. The  $\Omega$-calibrated  submanifolds 
(i.e.\ they satisfy $\cos\theta=1$) are minimal if $\Omega$ is closed (see \cite{[HL]}). 
In \cite{[LiSa]} we have got the following $\Omega$-integral isoperimetric inequality on 
a domain $D$ of $M$, that we display now in the most general case of $\Omega$, not necessarily closed, by following  the same proof,
\begin{equation}\label{isopintegral}
\left|\int_D\La{(}m\cos\theta\, \|H\|^2 -\langle \nabla^{\bot}H,\Phi\rangle\La{)}dV^*
-\int_D (\bar{\nabla}_H\Omega - H\lrcorner d\Omega)
\right|\leq \int_{\partial D}\sin\theta \|H\|dA^*.
\end{equation}
In this inequality, $dV^*$ and $dA^*$ are the volume and area forms for the induced metric 
$g_*=F^*\bar{g}$ on $M$,   $\bar{\nabla}$ is  and Levi Civita connection on $(\oM,\bar{g})$, and $\bar{\nabla}_H\Omega - H\lrcorner d\Omega$ is restricted to $M$ via pullback by $F$. 
The morphism  $\Phi:TM\to NM$,  with values on the normal bundle of $M$, and appearing in the inequality,
is defined by
\begin{equation}\label{Phi}
\bar{g}(\Phi(X),U)=\Omega(U,*X),
\end{equation}
for any $X\in T_pM$ and $U\in NM_p$, 
where $*:TM\to \wedge^{(m-1)}TM$ is the star operator.
If $\Omega$ is parallel, $\cos\theta>0$ on $\bar{D}$,  and $F$ has parallel mean curvature, from the 
above inequality we obtain
\begin{equation}\label{old2}
\|H\|\leq \frac{1}{m}\left(\frac{\sup_{\partial D}\sin \theta}{\inf_D\cos\theta}\right)
\frac{A^*(\partial D)}{V^*(D)}.
\end{equation}
In particular, if $\h(M,g_*)=0$ and $\cos\theta>\epsilon>0$ on $M$, then $H=0$. 
We may relax the later assumption on $\cos \theta$ giving conditions at infinity. 
Assuming $(M,g_*)$ is complete with nonnegative Ricci tensor, then for some constant 
$C_2>0$ we have  $\h(B_r,g_*) \leq C_2/r$, $\forall r>0$, where
 $B_r$ is the ball of radius $r$ centred at a fixed point of $(M,g_*)$ (cf.\ \cite{[ADR]}),
with corresponding distance function $r(x)$. 
If in addition we assume for some constant $\beta\in [0,1)$, 
 $\cos\theta (x) \geq C_1r^{-\beta}(x)$, these two conditions are sufficient
to obtain $\|H\|\leq C r^{\beta-1}(x)$ (see proof of Theorem 1.4 of \cite{[LiSa]}), and conclude $H=0$. Here $C$ and $C_1$ are positive constants.

A Riemannian product manifold $\bar{M}=M\times N$ has a natural  calibration defined by the 
volume of the projection onto the first component
$$\Omega((X_1,Y_1), \ldots, (X_m,Y_m))= Vol_M(X_1,\ldots, X_m).$$
In the case  of a graph submanifold, $\Gamma_f:M\to  M\times N$, $\Gamma_f(x)=(x,f(x))$, if
$(X^*_i,df(X^*_i))$ is a direct o.n. frame of $\Gamma_f$,  then
$$\cos\theta= \Omega((X^*_1, df(X^*_1)), \ldots, (X^*_m, df(X^*_m))
=\frac{1}{\sqrt{det_g(g + f^*h)}}>0,$$
where the determinant of the graph metric  $g_*=\Gamma_f^*\bar{g}=g+f^*h$ is taken
with respect to a 
diagonalizing $g$-o.n. frame of $f^*h$.  The condition   
$\cos\theta\geq \epsilon>0 $  is equivalent to $\|df\|^2$ to be bounded, and consequently
 the metrics $g_*$ and $g$ are equivalent. In this case,
$(M,g_*)$ has zero Cheeger constant if and only if  $(M,g)$ has so.
In (\cite{[LiSa]}) we obtain the conclusion $H=0$ for graph submanifolds with parallel mean curvature using the  isoperimetric inequality with the $\Omega$-angle, but we need  an extra condition on $\cos\theta$ at infinity as explained above. Hence, this approach applied to graphs seems weaker than the one in \cite{[S1]}, \cite{[S2]}, but the curvature conditions may be different.

In this work we generalize the  isoperimetric inequality (\ref{old}) for a graph submanifold $\Gamma_f$ in an ambient space a warped product of two Riemannian manifolds, 
$(M^m, g)$ and $(N^n, h)$, defined by a warping function $\rho:M\to (0, +\infty)$. We denote
this space by $(\tilde{M}, \tilde{g})$, where $\tilde{M}=M\times_{\rho} N$ is the product space 
$M\times N$ with the warped metric 
$$\tilde{g}=g+ \rho^2h,\quad\mbox{where}~\rho=e^{\psi}.$$ 
If $n=1$, some Bernstein type results for a graph  hypersurface in a warped product  have  
been obtained in \cite{[AD1],[AD2],[ADR],[Brendle],[Z1],[Z2], [Ra]}.
In this paper we  work in any codimension, and consider the Heinz estimation type problem
for the mean curvature.
This problem has been studied in \cite{[Im]} for $n=1$, in the context of weighted manifolds.
 We will extend Theorem 1 of \cite{[Im]} to higher codimensions,
using their concept of weighted Cheeger constant.

A graph submanifold is defined as an immersion $\Gamma_f:M\to \tilde{M}$,
$\Gamma_f(x)=(x,f(x))$, for a given map $f:M\to N$. It defines on $M$ a induced metric,
the graph metric
\begin{equation}\label{graphmetric}
 g_*(X,Y)=\Gamma_f^*\tilde{g}(X,Y)=g(x)(X,Y)+\tilde{h}(x)(df(X),df(Y)),
\end{equation}
where $\tilde{h}(x)=\rho^2(x)h(f(x))$ is a Riemannian metric on the pullback tangent 
bundle $f^{-1}TN$.
We denote by $df^*:f^{-1}TN\to TM$,  the adjoint morphism of $df:TM \to f^{-1}TN$ 
when we consider on $TM$ the graph metric $g_*$,  and on 
$f^{-1}TN$ the metric  $\tilde{h}$. 
The  endomorphisms, $Id-df df^*:f^{-1}TN\to f^{-1}TN$, and $Id-df^* df: TM\to TM$, symmetric
for the metrics $\tilde{h}$ and $g_*$, respectively,  are  both positive diffeomorphisms, 
with the same set of eigenvalues. These are useful diffeomorphisms.
For example, for any function $\Theta:M\to \mathbb{R}$, the following equality holds
\begin{equation}\label{gradientes}
(Id-df^*df)(\nabla^M \Theta)=\nabla^*\Theta,
\end{equation}
 where $\nabla^M\Theta$ and  $\nabla^*\Theta$  are the gradients with respect to 
$g$ and $g_*$, respectively.  We denote by $\bar{\nabla}df$ the Hessian of 
$f:(M,g)\to (N,h)$,
and consider the following vector fields $\Psi_*, W, W_1$ of $f^{-1}TN$ and $Z_1$ of $TM$,
given by 
\begin{eqnarray}
\Psi_* &=& df\left(\|df\|^2_{g_*}\,\nabla^M\psi +2\nabla^*\psi\right) \label{Psi}\\
W &=& trace_{g_*}\bar{\nabla}df \label{W}\\ 
Z_1 &=& df^*(W +\Psi_*) \label{Z1}\\
W_1 &=& (Id -dfdf^*)(W+\Psi_*). \label{W1}
\end{eqnarray}
\noindent
We  will prove in Lemma \ref{lemma 2.3} that
the mean curvature of $\Gamma_f$, $H=(H_M,H_N)$, is given by
\begin{equation} \label{meancurvature}
 mH=m(H_M,H_N)=(-Z_1,W_1)=(0, W+\Psi_*)^{\bot},
\end{equation}
where $(\cdot)^{\bot}$ denotes the orthogonal projection onto the normal bundle, 
and a minimal graph on $\tilde{M}$  is   defined by the equality 
\begin{equation}\label{minimal}
W_*:= W+\Psi_*=0.
\end{equation}

We state our first  main Theorem \ref{Main Theorem 1} for graph submanifolds
of $\tilde{M}$ 
that generalizes the case $\rho=1$ of \cite{[S2]}, and the case $n=1$ of \cite{[Im]}. 
We have the following relations between the angle of $\nabla^M\psi$ with the  
$M$ and $N$ components of the mean curvature 
\begin{equation}\label{twoAngles}
 g(H_M,\nabla^M\psi)=-\tilde{h}(H_N,df(\nabla^M\psi)).
\end{equation}
 Then we define the following set, that is empty if $f$ is constant along integral curves of $\nabla^M \psi$,
\begin{equation}\label{M^-}
M^-=\{x\in M: g(H_M,\nabla^M\psi)<0\}=\{x\in M: \tilde{h}(H_N,df(\nabla^M\psi))>0\}.
\end{equation}

The proof of Theorem \ref{Main Theorem 1} relies on  a key  formula 
that expresses $\|H\|^2$ as a weighted divergence of a suitable vector. We will show 
the following equality holds
\begin{equation}\label{KEYg}
e^{-\psi}div_g(e^{\psi}H_M)=div_g(H_M)+g(H_M,\nabla^M\psi)= -m\|H\|^2.
\end{equation}
The weighted Cheeger constant, introduced in \cite{[Im]}, is given by
\begin{equation}\label{heightedCheeger}
\h(M,g,\psi):=\inf_D\frac{A_{\psi}(\partial D)}{V_{\psi}(D)},
\end{equation}
where the ${\psi}$-weighted area and volume are considered on compact domains 
$\bar{D}=D\cup \partial D$ of
$(M,g)$ with smooth boundary, 
\begin{eqnarray}
A_{\psi}(\partial D)= \int_{\partial D}e^{\psi}dA, \quad&&\quad  
V_{\psi}(D)= \int_De^{\psi}dM.
\end{eqnarray}
 Clearly, if $\h(M,g)=0$, and  $\psi$ is bounded, then $\h(M,g,\psi)=0$. 
The ${\psi}$-weighted Cheeger constant has the following spectral property (see Section 3),
for $M$  compact with non-empty smooth boundary,
\begin{equation}\label{CheegerIneq}
\lambda_{\psi,1}(M)\geq \frac{1}{4}(\h(M,g,\psi))^2,
\end{equation}
where $\lambda_{\psi,1}(M)$ is the lowest eigenvalue of the
spectrum of the drift ${\psi}$-Laplacian,  $-\Delta_{\psi}u =-\Delta u- g(\nabla \psi, \nabla u)$,
with  Dirichlet boundary condition, $u=0$ on $\partial M$.
\begin{theorem}\label{Main Theorem 1}
 If $f:M\to N$ defines a graph submanifold $\Gamma_f:M\to M\times_{\rho}N$,
$\Gamma_f(x)=(x, f(x))$, 
with parallel mean curvature $H=(H_M,H_N)$,   then the following  estimates hold
\begin{equation}\label{weightedestimation}
m\|H\|\leq \h(M,g,\psi),
\end{equation}
\begin{equation}\label{estimationH}
 m\|H\| \leq  \h(M,g) +\sup_{M^-}\|\nabla^M\psi\|,
\end{equation} 
where $M^{-}$ is defined in (\ref{M^-}). In particular, if  $\h(M,g, \psi)=0$, or if 
$\h(M,g)=0$ and either $g(H_M,\nabla^M\psi)\geq 0$ for all $x\in M$ or
$\psi$ is bounded, then $H=0$.
\end{theorem}
\begin{corollary} \label{corollary} If $(M,g)$ is a closed Riemannian manifold then
any graph submanifold $\Gamma_f$ of $M\times_{\rho}N$ with parallel mean curvature is minimal.
\end{corollary}

We consider the particular case  $M$ is a complete non-compact  spherically symmetric Riemannian manifold, $M=M_{\tau}:=[0,+\infty)\times_{\tau}\mathbb{S}^{m-1}$ with metric 
$g=g_{\tau}:=dt^2+\tau^2(t)\sigma^2$, where $\tau\in C^2([0,+\infty))$ satisfies 
$\tau(0)=\tau''(0)=0$, $\tau'(0)=1$, $\tau(t)>0$ for all $t>0$, and $d\sigma^2$ is the Euclidean metric of the unit $(m-1)$-sphere. For $x=(t,\xi)$, let $r(x)$ be the distance function to the origin $o$ of $M_{\tau}$, the point that is identified with the class  
of all elements $(0,\xi)$, with $\xi \in \mathbb{S}^{m-1}$.
In next proposition we  build entire graph hypersurfaces with nonzero constant mean 
curvature in $M\times_{\rho}\mathbb{R}$ under suitable conditions on the warping functions 
$\tau$ and $\rho=e^{\psi}$, following a  similar construction as in \cite{[S1],[S2], [S3]}. Some details on spherically symmetric spaces with a density can be seen in \cite{[FS]}.
 Assume $\psi$ is a radial function on $M_{\tau}$, that is,  $\psi(x)=\Psi(r(x))$ where 
$\Psi:[0, +\infty)\to [0,+\infty)$ satisfies $\Psi'(0)=0$ (to make  $\psi$ of 
class $C^1$ at $x=o$). Consider the  function  $X:[0, +\infty)\to [0, +\infty)$ given by
$$X(t) =e^{\Psi(t)}(\tau(t))^{(m-1)}.$$
This function is positive for $t>0$ and  has a zero of finite order  $(m-1)$ at $t=0$.
 The function $\phi:[0, +\infty)\to [0, +\infty)$,
$$ \phi(t)=\frac{\int_0^tX(s)ds}{X(t)},$$
satisfies the equation $\phi' (t)= 1-\frac{X' (t)}{X(t)}\phi(t)$, with initial conditions
$\phi(0)=0$ and $\phi' (0)=1$.
\begin{proposition}\label{EXAMPLE} Assume the  infimum
$C_0:=\inf_{t\geq 0}\frac{1}{\phi(t)}$ is positive.
For each pair of constants $|c|<C_0$, and $d\in \mathbb{R}$,  
consider the functions, $\phi_c:[0,+\infty)\to (-1,1)$ defined by 
$ \phi_c(t) =c\phi(t)$, and $F=F_{c,d}:[0, +\infty)\to \mathbb{R}$  by
\begin{equation}\label{GRAPH}
F(t)= \LA{\int}_0^t \left(e^{-\Psi(s)}\frac{\phi_c(s)}{\sqrt{1-\phi_c(s)^2}}\right)ds +d.
\end{equation}
Then $f(x)=F(r(x))$ defines a graph hypersurface $\Gamma_f(x)=(x, f(x))$ that has
constant mean curvature $\frac{c}{m}$ in $\tilde{M}=M\times_{e^{\psi}}\mathbb{R}$.
 In particular,  $\h(M_{\tau}, g_{\tau}, \psi)\geq C_0$.
\end{proposition}
\noindent
If the infimum $C_0$ is taken at a point $t_0\in (0, +\infty)$, then $C_0$ is positive. 
Since $\lim_{t\to 0} \frac{1}{\phi(t)}=\lim_{t\to 0} \Psi'(t) +(m-1)\frac{\tau'}{\tau}(t)=+\infty$,  in order to have $C_0>0$ we only need to assume $\phi(t)$ 
bounded. This is the case when,  $M_{\tau}$ is the $m$-dimensional Hyperbolic space where 
$\tau(t)=\sinh(t)$ and $\psi=0$, giving $C_0=(m-1)=\h(\mathbb{H}^{m})$.\\[2mm]

The slices $M\times \{q\}$ are totally geodesic submanifolds of $\tilde{M}$, but
the $m$-form on $\tilde{M}$,
$$\Omega((X_1,U_1), \ldots, (X_m,U_m))=Vol_g(X_1, \ldots, X_m),$$
is not a parallel $m$-form if $\psi$ is not constant. However, we will show in Lemma 
\ref{closed not parallel} that $\Omega$ is a closed calibration  that calibrates the 
slices. Now we state our second main Theorem \ref{Main Theorem 2} that generalizes  
Theorem 1.4 of \cite{[LiSa]}, for graph submanifolds of $\tilde{M}=M\times_{\rho}N$ with
parallel mean curvature.
The key of the proof is the following divergence-type formula 
\begin{eqnarray}
  e^{-\psi}div_{g_*}(e^{\psi}\cos\theta H_M)&=&  div_{g_*}(\cos\theta H_M)
+g_*(\cos\theta H_M, \nabla^*\psi)\nonumber \\
&=& div_{g_*}(\cos\theta H_M))-(\tilde{\nabla}_H\Omega)(X^*_1, \ldots, X^*_m)\nonumber \\
&=& = -m\cos\theta \|H\|^2, \label{KEYcalibration}
\end{eqnarray}
where $X_i^*$ is a $g_*$-o.n. frame of $M$.
 We are assuming $M$ is oriented. If $(M,g)$ is complete, then so is $(M,g_*)$.
\begin{theorem}\label{Main Theorem 2} Assume $\Gamma_f$ has parallel mean curvature on 
$M\times_{\rho}N$, and $(M,g_*)$ is complete. Moreover, 
assume that for some constants $\beta\in [0,1)$ and $C_1>0$, 
 the $\Omega$-angle of $\Gamma_f$,
$$ \cos\theta:=\Omega(d\Gamma_f(X_1^*),\ldots, d\Gamma_f(X_m^*))= \frac{1}{\sqrt{det_g g_*}},$$ 
satisfies $\cos\theta\geq C_1 r^{-\beta}$, when $r\to +\infty$,  where
$r(x)$ denotes the distance function on $(M,g_*)$ from a fixed point $x_0$.
Furthermore, assume  for some nonnegative constants $\alpha\geq  \delta$,
 the $\psi$-Ricci tensor of $(M,g_*)$ is  bounded from below by $\alpha$,
and $\|\nabla^*\psi\|\leq \delta^{1/2}$. 
Then $\Gamma_f$ is a minimal submanifold.
\end{theorem} 
We will see in Corollary \ref{Corollary 3.1} of Section 3, under the above boundedness
conditions on  $\|\nabla^*\psi\|$ and on the
 $\psi$-Ricci tensor,   that $\h(M,g_*,\psi)=0$.

\section{Graphs in warped products}

We consider two Riemannian manifolds $(M^m, g)$ and $(N^n, h)$, and a function 
$\psi:M\to \mathbb{R}$, defining a Riemannian space  $(\tilde{M}, \tilde{g})$, 
where $\tilde{M}=M\times N$ is endowed with the warped metric 
$\tilde{g}=g+ e^{2\psi}h$. 
Let $\nabla^M$, $\nabla^N$ and $\tilde{\nabla}$ denote the Levi-Civita connections of 
$(M,g)$, $(N,h)$  and $(\tilde{M},\tilde{g})$ respectively.  
We recall the following properties of $\tilde{\nabla}$ (cf.\ \cite{[ONeill]}).
If $ X,Y$ are vector fields of $M$ and $U, W$ of $N$ then 
$$\tilde{\nabla}_XY=\nabla^M_XY,\quad\quad \tilde{\nabla}_XU=\tilde{\nabla}_UX=d\psi(X)U,$$
$$\tilde{\nabla}_UW=(\tilde{\nabla}_UW)^{Tan}+(\tilde{\nabla}_UW)^{Nor}=
\nabla_U^NW-\tilde{g}(U,W)\nabla^M\psi,$$
where for each $x\in M$ and $q\in N$, $Tan:T_{(x,q)}(M\times N)\to T_{(x,q)}(x\times N)$ and 
$Nor:T_{(x,q)}(M\times N)\to T_{(x,q)}(M\times q)$, are given by the projections,
$Tan(X,U)=U$, $Nor(X,U)=X$.\\[-2mm]

Given a map $f:M\to N$,  the graph submanifold, 
$\Gamma_f=\{(x,f(x)):x\in M\}\subset \tilde{M}$, 
can be seen as the embedding of $M$ by the graph map
$\Gamma_f:M\to \tilde{M}$, $\Gamma_f(x)=(x,f(x))$, that induces on $M$ the graph metric 
(\ref{graphmetric}), 
$g_*(X,Y)=g(X,Y)+\tilde{h}(x)(df(X),df(Y))$, where $\tilde{h}(x)=e^{2\psi(x)}h(f(x))$
is a Riemannian metric on the pullback tangent bundle $f^{-1}TN$.
We denote by $\nabla^{f^{-1}}$ the pullback connection of $\nabla^N$ by $f$
on $f^{-1}TN$.
\begin{lemma} \label{lemma 2.1} If $X,Y$ are smooth vector fields on $M$, and $U$ is a 
section of $f^{-1}TN$, then $(Y,U)$ is a section of $\Gamma_f^{-1}T\tilde{M}$ and
$$\tilde{\nabla}_X^{\Gamma^{-1}_f}(Y,U)=(\nabla^M_XY,\nabla_X^{f^{-1}}U)
+(-\tilde{h}(df(X),U)\nabla^M\psi \,,\, d\psi(X)U +d\psi(Y)df(X)).$$
\end{lemma}
\noindent
\em Proof. \em We take $U_i$ a local o.n.\ frame of $TN$, and write 
$U(x)=\sum_i\lambda_i(x,f(x))U_i(f(x))$, where $\lambda_i(x,q)=\lambda_i(x)$.
Hence, $(Y,U)=(Y,0)+ \sum_i\lambda_i(0,U_i)$, and we have
\begin{eqnarray*}
\lefteqn{\tilde{\nabla}_X^{\Gamma^{-1}_f}(Y,U)=
\tilde{\nabla}_{(X,df(X))}\La{(}(Y,0)+\sum_i \lambda_i(0,U_i)\La{)}}\\
&=& \tilde{\nabla}_{(X,df(X))}(Y,0)+
\sum_i\lambda_i\tilde{\nabla}_{(X,df(X))}(0,U_i)+d\lambda_i(X)(0,U_i)\\
&=&
\tilde{\nabla}_{(X,0)}(Y,0)+\tilde{\nabla}_{(0,df(X))}(Y,0)+
\sum_i \lambda_i\tilde{\nabla}_{(X,0)}(0,U_i)\\
&&+\sum_i\lambda_i\tilde{\nabla}_{(0,df(X))} (0,U_i)+ d\lambda_i(X)(0,U_i).
\end{eqnarray*}
Applying the above rules of the $\tilde{\nabla}$-connection, 
 we get the final expression.\qed\\[4mm]
Let $\nabla^*$ be the Levi Civita connection of $M$ for the  graph metric $g_*$.
The second fundamental form of $\Gamma_f$,  $\nabla^*d\Gamma_f:TM\times TM\to
N\Gamma_f$, takes values on the normal bundle  $N\Gamma_f\subset  \Gamma_f^{-1}T\tilde{M}$,
$${\nabla}^{*}d\Gamma_f(X,Y)= \tilde{\nabla}_X^{\Gamma_f^{-1}}(d\Gamma_f(Y))-d\Gamma_f(\nabla_X^*Y).$$
We denote by $(Y,U)^{\top}$ and $(Y,U)^{\bot}$ the $\tilde{g}$-orthogonal projections onto
$T\Gamma_f$ and  $N\Gamma_f$ respectively, and the Hessian of $f$ for the Levi-Civita connections $\nabla^M$ and $\nabla^N$ by
$$\bar{\nabla}df(X,Y)=\nabla^{f^{-1}}_X(df(Y))-df(\nabla^M_XY)
=\nabla^N_{df(X)}(df(Y))-df(\nabla^M_XY).$$
\begin{lemma} \label{lemma 2.2} 
$(1)$ $(Y,U)^{\bot}=(0,U-df(Y))^{\bot}$. 
 Hence, $(Y,U)^{\bot}=(0,df(Z))^{\bot}$ if and only if $(0,U)^{\bot}
=(0, df(Z+Y))^{\bot}$;
$(2)$ $(X,0)\in N\Gamma_f$ if and only if  $X=0$;
$(3)$ $(0,U)^{\bot}=0$ if and only if  $U=0$.
\end{lemma}
\noindent
\em Proof. \em (1)$(Y,U)^{\bot}=(Y,U)^{\bot}-(Y,df(Y))^{\bot}$. 
(2) If $(X,0)\in N\Gamma_f$, then $0=\tilde{g}((X,0),(X,df(X)))$
$=g(X,X)$. 
(3) If $(0,U)\in T\Gamma_f$ then $(0,U)=(X,df(X))$ for some $X$, and so $X=0$. \qed\\[1mm]

We define a $TM$-valued symmetric bilinear form,  $\Xi:T_xM\times T_xM \to T_xM$,
\begin{equation}\label{Xi}
\Xi(X,Y)=\tilde{h}(df(X),df(Y))\nabla^M\psi +g(\nabla^M\psi,X)Y+g(\nabla^M\psi,Y)X,
\end{equation}
and consider
 the adjoint linear morphism $df^*:f^{-1}TN\to TM$ of $df:TM\to f^{-1}TN$, considering  
$TM$ with the metric $g_*$ and $f^{-1}TN$ with $\tilde{h}$, that is, it is defined by 
$$g_*(df^*(U),X)=\tilde{h}(df(X),U).$$
We recall the vector fields $\Psi_*$, $W$, $Z_1$ and $W_1$ given in (\ref{Psi}), 
 (\ref{W}), (\ref{Z1}), and (\ref{W1}), and  set
\begin{eqnarray}\label{W*}
W_*&=& W+\Psi_*,\\
\Xi_*&=& trace_{g_*}\Xi.\label{Xi*}
\end{eqnarray}
We have
\begin{equation} \label{PsiXi*}
\Psi_*=df(\Xi_*), ~~\mbox{and}~Z_1=df^*(W_*).
\end{equation}
Let $X_i$  be a local  o.n. frame of $M$ with respect to $g$,
and set 
\begin{eqnarray*}
{g_*}_{ij} :=g_*(X_i,X_j)=  \delta_{ij} +\tilde{h}(x)
(df(X_i),df(X_j))
= \tilde{g}(d\Gamma_f(X_i),d\Gamma_f(X_j))=\tilde{g}_{ij}.
\end{eqnarray*}
Note that,
\begin{eqnarray}
 \sum_kg_*^{kr}\tilde{h}(W_*, df(X_k))&=& \sum_kg_*^{kr}g_*(df^*(W_*),X_k)
\nonumber \\
&=&  \sum_{ks}g_*^{kr} g(Z_1, X_s){g_*}_{sk}
=g(Z_1,X_r). \label{Z1Xr}
\end{eqnarray}

\begin{lemma}\label{lemma 2.3} We have
\begin{eqnarray}
{\nabla}^{*}d\Gamma_f(X,Y)
&=& (0, \bar{\nabla}df(X,Y) + df(\Xi(X,Y)))^{\bot}\label{hessGammaf}\\ 
mH &=& (0,W_*)^{\bot}=(-Z_1,W_1).\label{HZW}
\end{eqnarray}
In particular,
\begin{equation} \label{ineqHZ1}
m\|H\|\geq \|Z_1\|_g.
\end{equation}
Furthermore, $H=0$ if and only if $W=-\Psi_*$.
\end{lemma}
\noindent
\em Proof. \em  Applying Lemma \ref{lemma 2.1} to $d\Gamma_f(Y)=(Y,df(Y))$, and using the fact 
that the second fundamental form takes values on the normal bundle $N\Gamma_f$, we get 
\begin{eqnarray*}
{\nabla}^{*}d\Gamma_f(X,Y) &=& d\Gamma_f(\nabla^M_XY-\nabla^*_XY) +(0,\bar{\nabla}df(X,Y))\\
&& +\la{(}-\tilde{h}(df(X),df(Y))\nabla^M\psi~,~ d\psi(X)df(Y)+d\psi(Y)df(X)\la{)}\\
&=&(0,\bar{\nabla}df(X,Y))^{\bot}
+\la{(}-\tilde{h}(df(X),df(Y))\nabla^M\psi~,~ df(d\psi(X)Y\!+\!d\psi(Y)X)\la{)}^{\bot}\quad\quad\\
&=& (0, \bar{\nabla}df(X,Y) + df(\Xi(X,Y)))^{\bot},
\end{eqnarray*}
where in the last equality we used Lemma \ref{lemma 2.2}(1). Taking the $g_*$-trace 
of $ {\nabla}^{*}d\Gamma_f$, we get the first expression for the mean curvature
in (\ref{HZW}).
From (\ref{Z1}), (\ref{W1}), (\ref{W*}), (\ref{Xi*}) and (\ref{Z1Xr}),  we have
\begin{eqnarray*}
mH&=& (0,W_*)^{\bot}=(0,W_*)-(0,W_*)^{\top}\\
&=&(0,W_*)-\sum_{kr}g_*^{kr}\tilde{g}((0,W_*), (X_k,df(X_k)))(X_r,df(X_r))\\
&=&(0,W_*)-\sum_{kr}g_*^{kr}\tilde{h}(W_*, df(X_k))(X_r,df(X_r))\\
&=& (0, W_*)-\sum_rg(Z_1,X_r)(X_r,df(X_r))\\
&=& (0,W_*)-(Z_1,df(Z_1))= (-Z_1, W_1).
\end{eqnarray*}
Consequently,
$m^2\|H\|^2=\|Z_1\|^2_g+ \tilde{h}(W_1,W_1)\geq \|Z_1\|^2_g.$
Therefore, $H=0$ if and only if $Z_1=W_1=0$.  By  Lemma \ref{lemma 2.2}(3),
the later is equivalent to $W_*=0$.\qed\\[2mm]

We now choose $X_i$ that is a diagonalizing $g$-o.n.\ basis  of $\rho^2f^*h$
at a given point $x$,  that is,
\begin{equation}\label{eigen}
\rho^2(x)h(f(x))(df(X_i),df(X_j))=\lambda_i^2\delta_{ij},
\end{equation}
 with  $\lambda_1^2\geq \lambda_2^2\geq \ldots
\lambda_k^2\geq \lambda_{k+1}^2\geq\ldots\lambda_m^2\geq 0$, where $k$ is defined by
 $\lambda_k>0$ and $\lambda_{k+1}=0$ if $df(x)\neq 0$, otherwise we set $k=0$.
 Thus, $m-k$ is the dimension of the kernel of  $df$, and taking 
\begin{equation}\label{eigen*}
X^*_i:= \frac{X_i}{\sqrt{1+\lambda_i^2}}, \mbox{~~for~~} i=1,\ldots m,
\end{equation}
we obtain a $g_*$-o.n. frame. Hence,  $\{ X_1, \ldots, X_k\}$ span 
$df^*(f^{-1}TN)=(Kern~ df)^{\bot}$, and $\{X_{k+1},\ldots, X_m\}$ span 
$ Kern~df$.
The vector fields of $f^{-1}TN$,
\begin{equation}\label{eigenU}
U_i:= \lambda_i^{-1}df(X_i), \mbox{~~for~~} i\leq k,
\end{equation}
span $df(TM)$, a subspace which $\tilde{h}$-orthogonal complement is the kernel of $df^*$. 
We extend this $\tilde{h}$-o.n. system to give a $\tilde{h}$-o.n. basis of
$f^{-1}TN$, defining the subspace
$Kern~df^*=(df(TM))^{\bot} =span \{U_{k+1},\ldots, U_n\}$.
Recall that we are considering  $\lambda_{k+1}=\lambda_{k+2}=\ldots=\lambda_{\max\{m,n\}}=0$,
where $k\leq \min\{m, n\}$. Then we have,
 \begin{equation}\label{allEigen} 
\left\{\begin{array}{llll}
df(X^*_i)= \frac{\lambda_i}{\sqrt{1+\lambda_i^2}}U_i,&
df^* df(X_i)=\frac{\lambda_i^2}{{1+\lambda_i^2}}X_i,&
(Id-df^* df)(X_i)=\frac{1}{{1+\lambda_i^2}}X_i,&~i=1,\ldots, m,\\
df^*(U_i)=\frac{\lambda_i}{\sqrt{1+\lambda_i^2}}X^*_{i},&
 df df^*(U_i)=\frac{\lambda_i^2}{{1+\lambda_i^2}}U_i,&
 (Id- df df^*)(U_i)=\frac{1}{{1+\lambda_i^2}}U_i.&~i=1,\ldots, n.
\end{array}\right.
\end{equation}
It follows  that, $\|df\|^2_g=\sum_i\lambda_i^2$, and 
$ \|df\|^2_{*}= \|df^*\|^2_{*}= trace(dfdf^*)=trace (df^*df)
=\sum_{i}\frac{\lambda_i^2}{{1+\lambda_i^2}}$.
From the above identities (\ref{allEigen}) we  also derive the formula on the 
gradients (\ref{gradientes}).

Now we consider a morphism, 
$ Q_{\psi}:f^{-1}TN\to \mathbb{R}_M$,  where for $ U\in T_{f(x)}N$ it is given by
\begin{equation}  \label{Qpsi}
Q_{\psi}(U) = \tilde{h}\La{(}U\,,\, df(\nabla^*\psi)\La{)}.
\end{equation}
The morphism $Q_{\psi}$ is null  if $f$ is constant along the integral curves of 
$\nabla^*\psi$,
or of $\nabla^M\psi$, as a consequence of the identity (\ref{Qpsi*}) in  next 
Lemma \ref{BIGBANG}.   Using the
eigenvectors $X_i$, $X^*_i$, $U_i$ we  see that 
$
Q_{\psi}(\Psi_*)=\sum_j\frac{\lambda_j^2}{1+\lambda_j^2}
\left(\|df\|^2_*+\frac{2}{1+\lambda_j^2}\right)g(\nabla^M\psi,X_j)^2\geq 0$.

\begin{lemma}\label{BIGBANG} The following equalities hold
\begin{equation}\label{Qpsi*}
Q_{\psi}(U) =g(df^*(U),\nabla^M\psi) 
= \tilde{h}\La{(}(Id-df df^*)(U)\,,\, df(\nabla^M\psi)\La{)}.
\end{equation}
 Furthermore,
\begin{equation}\label{equalities}
Q_{\psi}(W_*) = g(Z_1, \nabla^M\psi) =\tilde{h}(W_1,df(\nabla^M\psi)).
\end{equation}
Hence, $M^-=\{x: Q_{\psi}(W_*)>0\}$, where $M^-$ is defined in (\ref{M^-}).
Then $M^-=\emptyset$ if and only if $Q_{\psi}(W)\leq -Q_{\psi}(\Psi_*)$ everywhere. 
\end{lemma}
\noindent 
\em Proof. \em The first equality of (\ref{Qpsi*}) follows from the identities
$$Q_{\psi}(U) =\tilde{h}\La{(}U\,,\, df(\nabla^*\psi)\La{)} =g_*(df^*(U),\nabla^*\psi)
=d\psi( df^*(U))=g(df^*(U),\nabla^M\psi).$$
Now we show the second.
\begin{eqnarray*}
\lefteqn{
\tilde{h}\La{(}(Id-df df^*)(U)\,,\, df(\nabla^M\psi)\La{)}= 
\sum_{i,j}\tilde{h}(U,U_i)\,d\psi(X_j)\,\tilde{h}\La{(}(Id-df df^*)(U_i)\,,\, df(X_j)\La{)}}\\
&=& \sum_{i\leq n,~j\leq k}\tilde{h}(U,U_i)\,d\psi(X_j)\,\frac{\lambda_j}{1+\lambda_i^2}\tilde{h}\La{(}U_i\,,\, U_j)
= \sum_{i\leq k} \tilde{h}(U,U_i)\,d\psi(X_i)\,\frac{\lambda_i}{1+\lambda_i^2}\\
&=& \sum_i \tilde{h}(U,df(X^*_i))\,d\psi(X^*_i)
 = \tilde{h}(U,df(\nabla^*\psi)) =Q_{\psi}(U).
\end{eqnarray*}
Then, $Q_{\psi}(W_*) =g_*(df^*(W_*),\nabla^*\psi)
=g_*(Z_1,\nabla^*\psi)=g(Z_1, \nabla^M\psi)$.
The proof of last equality in (\ref{equalities}) is a direct consequence of $H$ to take values on the normal bundle and Lemma \ref{lemma 2.3}, 
and so $\tilde{g}((-Z_1,W_1), (\nabla^M\psi,df(\nabla^M\psi)))=0$.\qed

If $X_i^*$ is a $g_*$-o.n. basis of $T_xM$, then
\[
\langle d\Gamma_f(\cdot),\tilde{\nabla}^{\Gamma_f^{-1}} H\rangle =
\sum_i\tilde{h}(d\Gamma_f(X_i^*), \tilde{\nabla}_{X_i^*}^{\Gamma_f^{-1}} H )
=-\sum_i\tilde{h}(\nabla^*_{X_i^*}d\Gamma_f(X_i^*), H)\]
Hence,
\begin{equation}\label{key}
\langle d\Gamma_f(\cdot),\tilde{\nabla} H\rangle = -m \|H\|^2.
\end{equation}

Now we prove main Theorem \ref{Main Theorem 1}. We denote by  $\nabla^{\bot}$ the  
covariant derivative on $N\Gamma_f$.\\[4mm]
\em Proof of  ~{\bf Theorem \ref{Main Theorem 1}}. \em 
From equalities (\ref{W*}), (\ref{Xi*}) and (\ref{PsiXi*}),
 $W_1=W+df(\Xi_*-Z_1)$, and  from Lemma \ref{lemma 2.3}, 
$mH= (-Z_1,W_1)$. Thus,  by Lemma \ref{lemma 2.1},
\begin{eqnarray*}
\lefteqn{m \nabla_{X}^{\Gamma_f^{-1}}H=}\\
&=& (-\nabla^M_XZ_1 ~,~ \nabla_{X}^{f^{-1}}W_1)
 +\La{(}-\tilde{h}(df(X),W_1)\nabla^M\psi~ , ~d\psi(X)W_1+d\psi(-Z_1)df(X)\La{)}\\
&=& -(\nabla_X^MZ_1, df(\nabla_X^MZ_1)) +\LA{(}0~,~ \nabla_X^{f^{-1}}W+\bar{\nabla}_Xdf(\Xi_*-Z_1)
+df(\nabla_X^M\Xi_*)\LA{)}\\
&& +\LA{(}-\tilde{h}(df(X),W_1)\nabla^M\psi~ , ~d\psi(X)W_1+d\psi(-Z_1)df(X)\LA{)}.
\end{eqnarray*}
By Lemma \ref{lemma 2.2} (1), if $\nabla^{\bot}H=0$ then 
\begin{eqnarray*}
&&\LA{(}0~,~ \nabla_X^{f^{-1}}W+\bar{\nabla}_Xdf(\Xi_*-Z_1)
+df(\nabla_X^M\Xi_*)\LA{)}^{\bot}\\
&& +\LA{(}0~,~ d\psi(X)W_1+d\psi(-Z_1)df(X)  +\tilde{h}(df(X),W_1)df(\nabla^M\psi)\LA{)}^{\bot}=0,
\end{eqnarray*}
and so,  by Lemma \ref{lemma 2.2}(3)
\begin{eqnarray*}
\lefteqn{ \nabla_X^{f^{-1}}W+\bar{\nabla}_Xdf(\Xi_*-Z_1) +df(\nabla_X^M\Xi_*)=}\\
&=&  -d\psi(X)W_1+d\psi(Z_1)df(X) -\tilde{h}(df(X),W_1)df(\nabla^M\psi).
\end{eqnarray*}
Thus,
\begin{eqnarray*}
m \nabla_{X}^{\Gamma_f^{-1}}H&=&
-(\nabla_X^MZ_1, df(\nabla_X^MZ_1))
 +\La{(}-\tilde{h}(df(X),W_1)\nabla^M\psi~ , 
-\tilde{h}(df(X),W_1)df(\nabla^M\psi)\La{)}\\
&=& d\Gamma_f\La{(}-\nabla_X^MZ_1-\tilde{h}(df(X),W_1)\nabla^M\psi\La{)}.
\end{eqnarray*}
Now we fix a point $x_0$ and take $X_i$ s.t. $\nabla^MX_i(x_0)=0$. Then, by 
(\ref{Z1Xr}) at $x_0$
$$
\nabla^M_{X_i}Z_1=\sum_{kr}d\left( g_*^{kr}\tilde{h}(W_*, df(X_k))\right)\!(X_i)\,X_r,$$
and so, for each $i,j$
\begin{eqnarray*}
\tilde{g}\LA{(}(\nabla_{X_i}^{M}Z_1,df(\nabla_{X_i}^{M}Z_1))\,,\, (X_j,df(X_j)\LA{)}&=&
\sum_{kr} d\left( g_*^{kr}\tilde{h}(W_*, df(X_k)\right)\!(X_i)\tilde{g}_{rj}.
\end{eqnarray*}
In particular, at $x_0$
\begin{eqnarray*}
div_{g}(Z_1) &=& \sum_rg(\nabla^M_{X_r}Z_1, X_r)=
\sum_{kr} d\La{(}g_*^{kr}\tilde{h}(W_*, df(X_k))\La{)}(X_r).
\end{eqnarray*}
Therefore, at $x_0$, 
\begin{eqnarray*}
\lefteqn{-m^2\|H\|^2 = m\langle \tilde{\nabla}^{\Gamma_f^{-1}}H, d\Gamma_f\rangle =
\sum_{ij}{g}_{*}^{ij}\tilde{g}(\tilde{\nabla}^{\Gamma_f^{-1}}_{X_i}H, d\Gamma_f(X_j))}\\
&=&
\sum_{ij}{g}_{*}^{ij}\tilde{g}\left(-d\Gamma_f\La{(}\nabla_{X_i}^MZ_1+\tilde{h}(df(X_i),W_1)\nabla^M\psi\La{)},  d\Gamma_f(X_j)\right)\\
&=&\sum_{ij}{g}_{*}^{ij}\tilde{g}\left( -\nabla_{X_i}^MZ_1, df(-\nabla_{X_i}^MZ_1)),
(X_j,df(X_j))\right)\\
&& + \sum_{ij}{g}_{*}^{ij}\left( -\tilde{h}(df(X_i),W_1)\tilde{g}
(d\Gamma_f(\nabla^M\psi),d\Gamma_f(X_j))\right)\\
&=& -\sum_{ijkr}{g}_{*}^{ij} d\left( g_*^{kr}\tilde{h}(W_*, df(X_k)\right)\!(X_i)\tilde{g}_{rj}\\
&& -\sum_{ij} g_*^{ij}\tilde{g}(d\Gamma_f(df^*(W_1)),d\Gamma_f(X_i))\tilde{g}(d\Gamma_f(\nabla^M\psi),
d\Gamma_f(X_j))\\
&=& -div_g(Z_1) - \tilde{g}(d\Gamma_f(df^*(W_1)), d\Gamma_f(\nabla^M\psi))\\
&=& -div_g(Z_1)-g_*(df^*(W_1),\nabla^M\psi)
 = -div_g(Z_1) -\tilde{h}(W_1, df(\nabla^M\psi)).
\end{eqnarray*}
Hence, by Lemma \ref{BIGBANG},
\begin{equation}\label{equationH}
m^2\|H\|^2 = div_{g}(Z_1) + g(Z_1, \nabla^M\psi) = e^{-\psi}div_g(e^{\psi}Z_1),
\end{equation}
where $\|H\|$ is constant by assumption on parallel mean curvature.
Weighted integration over $D$ of $e^{-\psi}div_g(e^{\psi}Z_1)$, using Stokes', Schwartz inequality  and applying 
Lemma \ref{lemma 2.3}, gives
\begin{eqnarray*}
m^2\|H\|^2V_{\psi}(D)&=& \int_D e^{-\psi}div_g(e^{\psi}Z_1)e^{\psi}dM 
=\int_D div_g(e^{\psi}Z_1)dM\\
&=&
\int_{\partial D}e^{\psi}g(Z_1,\nu)dA
\leq \int_{\partial D} \|Z_1\|_{g}e^{\psi}dA
\leq m\|H\|A_{\psi}(\partial D),
\end{eqnarray*}
where $\nu$ is the unit outward of $(\partial D,g)$. Hence
$$m\|H\|\leq \frac{A_{\psi}(\partial D)}{V_{\psi}(D)},$$
and  (\ref{weightedestimation}) of the Theorem is proved. 
Usual integration of $ div_{g}(Z_1) + g(Z_1, \nabla^M\psi)$, gives
\begin{eqnarray*}
m^2\|H\|^2V(D) &\leq& m\|H\|A(\partial D) +\int_{D\cap M^-}g(Z_1, \nabla^M\psi)dM\\
&\leq & m\|H\| A(\partial D) +m\|H\|\,sup_{M^-}\|\nabla^M\psi\|V(D).
\end{eqnarray*}
Hence $m\|H\|\leq \frac{A(\partial D)}{V(D)}+ \,sup_{M^-}\|\nabla^M\psi\|$, and  (\ref{estimationH}) is proved.\qed\\[4mm]

\noindent
\em Proof of  ~{\bf Proposition \ref{EXAMPLE} }. \em  
For any Riemannian manifold $(M,g)$, if  $N=\mathbb{R}$  and
$f:M\to \mathbb{R}$,  the unit normal of  $\Gamma_f$ in $\tilde{M}=M\times_{\rho}
\mathbb{R}$  is given by
$$ N=\frac{(-\nabla f, e^{-2\psi})}{\sqrt{e^{-2\psi}+ \|\nabla f\|^2_g}},$$
where $\nabla f=\nabla^M f$, and the mean curvature $H$ satisfies the equation
$$m\, \tilde{g}(H,N)=div_{-\psi}\left(\frac{\nabla f}{\sqrt{e^{-2\psi}+\|\nabla f\|^2}}
\right) :=e^{-\psi}div_g\left(e^{\psi}\frac{\nabla f}{\sqrt{e^{-2\psi}+\|\nabla f\|^2}}
\right).$$
Hence $m\,\tilde{g}(H,N)=c$ if and only if
$$ div_g\left(\frac{e^{\psi}\nabla f}{\sqrt{e^{-2\psi}+\|\nabla f\|^2}}
\right)=c\, e^{\psi}.$$
Now we are assuming  $(M,g)=(M_{\tau}, g_{\tau})$, and so 
$r(t,\xi)=t$, and $\nabla r(t,\xi)=\partial_t(t,\xi)$. Moreover,  
for any radial function $A(r)$ we have  
$$div_g(A(r)\nabla^M r)=A'(r)  + A(r)(m-1)(\log \tau)'=
A'(t)+ A(t)(m-1)\frac{\tau' (t)}{\tau(t)}.$$
We have  $f(x)=F(r(x))$, $\psi(x)=\Psi(r(x))$,
and $X(t) = e^{\Psi(t)}\,\tau(t)^{m-1}$. 
Then we  consider the following function
$$ \xi(t) = \frac{ F'(t) }{\sqrt{e^{-2\Psi(t)} + (F'(t))^2}}\in (-1, 1). $$
It follows that,  $m\tilde{g}(H,N)=c$
if and only if 
$ c=\xi' +  ( \Psi' + (m-1)(\ln\tau)' )\xi$, that is
$$ \xi' = c - (\ln X)' \xi.$$
Solutions $\xi$ of this equation are of the form 
$$ \xi=  \phi_c (t) =c\,\frac{\int_0^t X(s)ds}{X(t)}.$$
Since this function must take values in  $ (-1,1)$, we conclude 
$|c|<C_0$, and if $c\neq 0$,  $c\phi_c(t)>0$ for $t>0$. Furthermore, 
$\phi_c(0)=c\, \frac{d^{m-1}X}{dt^{m-1}}(0)/ \frac{d^{m}X}{dt^{m}}(0)=0$
and $\phi'_c(0)=c$.\qed\\[4mm]

Theorem 1.2 is a generalization of Theorem 1.4 of \cite{[LiSa]}, so we will give a 
sketch of the proof,  detailed in general, except on some statements that are  easy 
to follow from the references that will be indicated along the proof. First we prove 
some properties of the calibration $\Omega$.
\begin{lemma}\label{closed not parallel} The $m$-form  $\Omega$ is a
  closed  calibration of $\tilde{M}$ that calibrates the slices. It is
parallel if and only if $\psi$ is constant.
Moreover, if $X_i^*$ is a $g_*$-o.n. frame of $M$, we have 
\begin{eqnarray}\label{DOmega}
(\tilde{\nabla}_{H}\Omega)(d\Gamma_f(X_1^*),\ldots, d\Gamma_f(X_m^*))&=& 
\cos\theta\tilde{h}(H_N,df(\nabla^M\psi)) \nonumber\\
&=& -\cos\theta g(H_M,\nabla^M\psi)=-
\cos\theta g_*(H_M, \nabla^*\psi).\quad\quad
\end{eqnarray}
\end{lemma}
\noindent \em Proof. \em 
To see that $\Omega$ is a  closed $m$-form  we use a  $\tilde{g}$-o.n frame of the form $(X_i,0),(0,W_{\alpha})$ on $\tilde{M}$, being $X_i$ a direct o.n. frame of $(M,g)$, with 
$\nabla^MX_i(x_0)=0$ at a given point $x_0
\in M$.
Using  the properties of $\tilde{\nabla}$, we see that the covariant derivatives of $\Omega$ 
in the directions of all these vector fields vanish except for,
\begin{eqnarray*}
(\tilde{\nabla}_{(0,W_{\alpha})}\Omega) ((X_1,0),\ldots, (X_{i-1},0),(0,W_{\beta}),
(X_{i+1},0), \ldots, (X_m,0))\\
=
d\psi(X_i)\,\tilde{h}(W_{\alpha}, W_{\beta})
=\delta_{\alpha \beta}\,g(\nabla^M\psi,X_i).
\end{eqnarray*}
Thus, $\Omega$ is parallel only if $\psi$ is constant. To conclude $d\Omega=0$ we only need 
to check the following components of $d\Omega$,  
\begin{eqnarray*}
\lefteqn{
\pm d\Omega((X_1,0), \ldots, (X_{m-1},0),(0,W_{\alpha}),(0, W_{\beta}))=}\\
&=&
(\tilde{\nabla}_{(0,W_{\alpha})}\Omega) ((X_1,0),\ldots, (X_{m-1},0),(0,W_{\beta}))
-(\tilde{\nabla}_{(0,W_{\beta})}\Omega) ((X_1,0),\ldots, (X_{m-1},0),(0,W_{\alpha}))\\
&=&(d\psi(X_m)-d\psi(X_m))\delta_{\alpha \beta}=0.
\end{eqnarray*}
Note that $\tilde{\nabla}_H\Omega(d\Gamma_f(X^*_1),\ldots, d\Gamma_f(X^*_m))$ is 
independent of the
$g_*$-o.n. basis $X^*_i$ we take.
 We may assume that $X^*_i$ is defined by (\ref{eigen*}), where $X_i$ is a diagonalizing
$g$-o.n. basis of $T_{x_0}M$ of $\rho^2f^*h$, and  we chose
 $W_{\alpha}$ a local frame of $N$ that at $f(x_0)$ is given by $U_{\alpha}$ defined 
 as in (\ref{eigenU}) with respect to $X_i$.  Then,
 writing $H=(H_M,H_N)=\sum_i g(H_M,X_i)X_i+ 
\sum_{\alpha}\tilde{h}(H_N,U_{\alpha})U_{\alpha}$, we see that
\begin{eqnarray*}
\lefteqn{
(\tilde{\nabla}_{H}\Omega)(d\Gamma_f(X^*_1),\ldots, d\Gamma_f(X^*_m)) =}\\
&=&(\tilde{\nabla}_{H}\Omega)((X_1, df(X_1)),\ldots, (X_m,df(X_m)))\frac{1}{\sqrt{\Pi_i(1+\lambda_i^2)}}\\
&=& \sum_{\alpha,i}\tilde{h}(H_N,U_{\alpha})
(\tilde{\nabla}_{U_{\alpha}}\Omega)((X_1,0 ),\ldots, (X_{i-1},0), (0,df(X_i)), (X_{i+1},0)
\ldots, (X_m,0) )\cos\theta\\
&=& \sum_{i,\alpha}\tilde{h}(H_N,U_{\alpha})\tilde{h}(df(X_{i}), U_{\alpha})d\psi(X_{i})\cos\theta = \tilde{h}(\cos\theta H_N, df(\nabla^M\psi)).
\end{eqnarray*}
Lemma \ref{BIGBANG} proves the other equalities in (\ref{DOmega}).\qed\\[4mm]
\noindent
\em Proof of ~{\bf Theorem \ref{Main Theorem 2}. } \em 
Following the proofs of Theorems 1.2 and  1.3 of  \cite{[LiSa]}, we consider
the morphism $\Phi:TM\to N\Gamma_f$, defined in (\ref{Phi}). By Lemma 2.1 of \cite{[LiSa]},
$\|\Phi(X)\|^2\leq \sin^2\theta \|X\|^2_{g_*}$. We  define a vector field $Z$ on $M$ by the equality, 
$$g_*(Z,X)=\tilde{g}(\Phi(X),H).$$
Then, $\|Z\|_{g_*}\leq \sin \theta \|H\|$. Moreover, since $*X^*_i= (-1)^{i-1}X_1^* \wedge
X_{i-1}^*\wedge X_{i+1}^*\wedge \ldots X_m^*$, then
\begin{eqnarray*}
g_*(Z,X_i^*) &=& (-1)^{i-1}\Omega(H,d\Gamma(X_1^*),
d\Gamma_f(X_{i-1}^*), d\Gamma_f(X_{i+1}^*), \ldots d\Gamma_f(X_m^*))\\
&=&\frac{1}{\Pi_{s\neq i} \sqrt{1+\lambda_s^2}}\Omega((X_1,0)), \ldots , (H_M,0), \ldots , (X_m,0))\\
&=& g(H_M,X_i){\cos\theta}{\sqrt{1+\lambda_i^2}}\,
\Omega((X_1,0), \ldots, (X_m,0))
=g(\cos\theta H_M,X_i){\sqrt{1+\lambda_i^2}},
\end{eqnarray*}
that is 
$$Z= \sum_i g_*(Z,X^*_i)X^*_i= g(\cos\theta H_M,X_i)X_i=\cos\theta H_M.$$
 Furthermore, as in  Theorems 1.2 and  1.3 of  \cite{[LiSa]},
$$ div_{g_*}(Z)=-\tilde{g}(\delta\Phi, H)+\sum_i\tilde{g}(\Phi(X^*_i),\tilde{\nabla}^{\bot}_{
X^*_i}H).$$
If $N_{\alpha}$ is an
$\tilde{g}$-o.n. frame of $N\Gamma_f$, and identifying $X^*_i$  with $d\Gamma_f(X^*_i)$  then, as in the proof of Lemma 2.2 of \cite{[LiSa]},
we have
$$\delta\Phi =-\tilde{\nabla}_{X_i^*}\Phi(X_i^*)=
 \sum_{\alpha} d\Omega(N_{\alpha},X_1^*, \ldots, X_m^*)N_{\alpha}
-(\tilde{\nabla}_{N_{\alpha}}\Omega)(X_1^*, \ldots,X_m^*)N_{\alpha}+m\cos\theta\, H.$$
Hence, for $\Gamma_f$ with parallel mean curvature, 
$$
 div_{g_*}(Z) =-m\cos\theta\| H\|^2+\tilde{\nabla}_{H}\Omega(X_1^*, \ldots,X_m^*).$$
Therefore, by Lemma \ref{closed not parallel}
$$div_{g_*}(\cos\theta H_M)= -m\cos\theta \|H\|^2-g_*(\cos\theta H_M,\nabla^*\psi),$$ 
and we have proved equation (\ref{KEYcalibration}).
Weighted integration of  (\ref{KEYcalibration})  on a compact domain $D$ of $(M,g_*)$,   gives
$$ m\|H\|^2\int_D \cos\theta e^{\psi}dM^*=
-\int_{\partial D}g_*(\cos\theta H_M,\nu)e^{\psi}dA^*,$$
where $dM^*$ is the volume element of $(M,g_*)$, and $dA^*$ of $\partial D$, for the induced metric,
and $\nu$ the outward unit normal to $\partial D$.
Thus,
$$ m\|H\|^2(\inf_D\cos\theta)\, V_{\psi}^*(D)\leq 
 \|H\| A_{\psi}^*(\partial D).$$
If we take $D$ a domain of geodesic ball $B_r$ of $(M,g_*)$ of radius $r$, with 
$\bar{D}\subset B_r$, we have obtained
$$ m\|H\|\inf_{B_r}\cos\theta \leq  
\frac{A_{\psi}^*(\partial D)}{V_{\psi}^*(D)}.$$
Hence,
$$m\|H\|\inf_{B_r}\cos\theta\leq \h(B_r,g_*,\psi).$$
Using the assumption on $\cos\theta$, we have
$$m \|H\|\leq  C_1^{-1}r^{\beta}\h(B_r,g_*,\psi).$$
Under  the boundedness conditions on $\mathrm{Ricci}_{\psi,g_*}$ and 
$\|\nabla^*\psi\|_{g_*}$ we have
from  Corollary \ref{Corollary 3.1} of Section 3, $\h(B_r,g_*,\psi)\leq C'/r$. Hence
$$ \|H\|\leq  C r^{\beta-1}.$$
Making $r\to +\infty$ we get $H=0$
.\qed\\[2mm]

\section{Weighted Cheeger inequality}
We  prove the weighted Cheeger inequality (\ref{CheegerIneq})
that generalizes the well known Cheeger inequality (cf.\ \cite{[Cha]},  Section IV, Theorem 3).
Recall that the eigenvalue problem of the $\psi$-Laplacian,  
$-\Delta_{\psi}u=-\Delta u - g(\nabla u, \nabla \psi)$, on a smooth compact Riemannian
manifold $M$ with boundary, with Dirichlet boundary condition, $u=0$ on $\partial M$, 
 consists on a discrete sequence 
$$ 0<\lambda_{\psi,1}<\lambda_{\psi, 2}\leq \lambda_{\psi, 3}\leq \ldots \to +\infty$$
and each eigenvalue has a variational characterization of Rayleight type (cf.\ 
\cite{[LuRow]}). 
We have for the principal eigenvalue $\lambda_{\psi,1}$, 
$$ \lambda_{\psi, 1}(M)= \inf_{u\in C_0^{\infty}(M)}
\frac{\int_M\|\nabla u\|^2 e^{\psi}dM}{\int_M u^2 e^{\psi}dM},$$
and the infimum  is attained at a principal eigenfunction $u$.
\begin{theorem} \label{Theorem 3.1}
$$\lambda_{\psi,1}(M)\geq \frac{1}{4}(\h(M,g,\psi))^2.$$
\end{theorem}
\noindent
\em Proof. \em 
We follow the proof of Theorem 3 of \cite{[Cha]}, Section IV, for the case $\psi=0$.
Let  $u$ be a principal eigenfunction, normalized s.t. $\int_Mu^2e^{\psi}dM=1$. 
Then  $u>0$ on $M$ and vanishes on $\partial M$. Set 
$M(t)=\{x\in M: u^2(x)>t\}$, $\Sigma(t)=\{x\in M: u^2(x)=t\}$, and
$$V_{\psi}(t)=V_{\psi}(M(t))=\int_{M(t)}e^{\psi}dM,\quad\quad 
A_{\psi}(t)=A_{\psi}(\Sigma(t))=\int_{\Sigma(t)}e^{\psi}dA.$$ 
We first recall the co-area formula (cf.\
 \cite{[Cha]}, Sec IV., Theorem 1) for $f=u^2$ and any function $h\in L^1(M)$,
$$ \int_M h\|\nabla u^2\|dM=\int_0^{+\infty}ds\int_{\Sigma(s)}hdA.$$ 
Considering $h=e^{\psi}\chi_{M(t)}\|\nabla u^2\|^{-1}$,
where $\chi$ is the characteristic function,  we get
$$V_{\psi}(t)=\int_{M}\chi_{M(t)}e^{\psi}dM=
\int_0^{+\infty}\!\!\!ds\int_{\Sigma(s)}\chi_{M(t)}\|\nabla u^2\|^{-1}e^{\psi}dA=
\int_t^{+\infty}\!\!\!ds\int_{\Sigma(s)}\|\nabla u^2\|^{-1}e^{\psi}dA
$$
Hence $ V'_{\psi}(t)= -\int_{\Sigma(t)}\|\nabla u^2\|^{-1}e^{\psi}dA$.  
Again, by the co-area formula with $h=e^{\psi}$,
$$ \int_M\|\nabla u^2\|e^{\psi}dM=\int_0^{+\infty} A_{\psi}(s)ds. $$
Since $L^2_{\psi}(M)$ is an Hilbert space, we have
\begin{eqnarray*}
\left(\int_M\|\nabla u^2\|e^{\psi}dM\right)^2&=& 4\left(\int_Mu\|\nabla u\|e^{\psi}dM\right)^2\\
&\leq& 4\left(\int_M\|\nabla u\|^2e^{\psi}dM\right) \left(\int_M u^2e^{\psi}dM\right)= 4
\int_M\|\nabla u\|^2e^{\psi}dM.
\end{eqnarray*}
Hence,
$$\lambda_{\psi,1}=\int_M\|\nabla u\|^2e^{\psi}dM\geq 
\frac{1}{4}\left(\int_M\|\nabla u^2\|e^{\psi}dM\right)^2$$
Now, using that $tV_{\psi}(t)$ vanish at $t=0$ and $t=+\infty$, and the co-area formula 
with  $h=u^2\|\nabla u^2\|^{-1}e^{\psi}$, we have 
\begin{eqnarray*}
\int_M\|\nabla u^2\|e^{\psi}dM &=& \int_0^{+\infty} A_{\psi}(s)ds\geq 
\int_0^{+\infty} \h(M(s),g,\psi)V_{\psi}(s)ds\\
&\geq& \h(M,g,\psi)\int_0^{+\infty}V_{\psi}(s)ds= 
-\h(M,g,\psi)\int_0^{+\infty}sV'_{\psi}(s)ds
\\
&=&\h(M,g,\psi)\int_0^{+\infty}s ds\int_{\Sigma(s)}\|\nabla u^2\|^{-1}e^{\psi}dA\\
&=& \h(M,g,\psi)\int_Mu^2e^{\psi}dM=\h(M,g,\psi).
\end{eqnarray*}
Consequently $ \lambda_{\psi,1}\geq \frac{1}{4}(\h(M,g,\psi)^2$, 
and the theorem is proved. \qed\\[4mm]
On a weighted manifold $(M,g, e^{\psi})$, it is defined the $\psi$-Ricci tensor, 
$\mathrm{Ricci}_{\psi,g}(X,Y):=\mathrm{Ricci}_g -Hess(\psi)$.
We  recall the following comparison result due to Setti (\cite{[Setti]}, Theorem 4.2).
In our notation $\mathrm{Ricci}_{\psi,g}=S_{\omega}$ given in (2.2) of \cite{[Setti]},
where $\omega=e^{\psi}$,  and $R_{\omega}$ defined in \cite{[Setti]} corresponds to  
$\mathrm{Ricci}_{\psi,g}-d\psi
\otimes d\psi$. The  conditions on $\mathrm{Ricci}_{\psi,g}$ and $\|\nabla\psi\|$ 
stated in the next theorem implies $R_{\omega}\geq (\alpha-\delta)$. 
\begin{theorem} If $(M^m,g,e^{\psi})$ is a complete weighted manifold, and
$\mathrm{Ricci}_{\psi,g}\geq \alpha$, and $\|\nabla \psi\|\leq \delta^{1/2}$, where 
$\alpha,\delta$ are constants, then, the first eigenvalue for the $\psi$-Laplacian on  a geodesic ball of radius $r$ on $M$, $B_r$, satisfies
$$\lambda_{\psi,1}(B_r)\leq \lambda_1(B^0_r),$$
where $\lambda_1(B^0_r)$ is the first eigenvalue for the 0-Laplacian  of a geodesic ball of 
radius $r$ with Dirichlet boundary condition on a $(m+1)$-dimensional space form of 
constant sectional curvature $(\alpha-\delta)/m$.
\end{theorem}
\noindent 
Obviously we are assuming $\delta\geq 0$. 
If $\alpha -\delta\geq 0$ then for some constant $C'>0 $, 
$\lambda_1(B^0_r) \leq \frac{C'}{r^2}$
(cf.\ \cite{[ADR]} or proof of Proposition 4.2 in \cite{[LiSa]}), as a consequence of 
Cheng's comparison result \cite{[Cheng]}). Therefore, we have the following conclusion.
\begin{corollary}\label{Corollary 3.1} In the previous theorem,  
if $\alpha\geq \delta$, then, for each ball of radius $r$, we have
$\h(B_r,g,\psi)\leq C/r$. In particular,   $\h(M,g,\psi)=0$.
\end{corollary}
\noindent
\em Proof. \em   By  previous  theorem
 $\lambda_{\psi,1}(B_r)\leq \lambda_1(B^0_r)\leq C' /r^2$, 
Then,  applying Theorem \ref{Theorem 3.1},
 $\h(M,g,\psi)\leq \h(B_r,p,\psi)\leq C/r\to 0 $, when $r\to +\infty$.

\end{document}